\documentclass[11pt]{article}
\usepackage{amssymb,latexsym, amsmath}
\usepackage{amscd, amsthm}

\newcommand{\be}{\begin{eqnarray}}
\newcommand{\ee}{\end{eqnarray}}
\newcommand{\bew}{\begin{eqnarray*}}
\newcommand{\eew}{\end{eqnarray*}}
\makeatletter \@addtoreset{figure}{section}
\def\thefigure{\thesection.\@arabic\c@figure}
\def\fps@figure{h, t}
\@addtoreset{table}{bsection}
\def\thetable{\thesection.\@arabic\c@table}
\def\fps@table{h, t}
\@addtoreset{equation}{section}

\newif\ifamsfonts
\amsfontstrue \ifamsfonts
\font\twlbbb=msbm10 scaled\magstep1 \font\egtbbb=msbm8
\font\sixbbb=msbm6
\newfam\mathbbfam
\textfont\mathbbfam=\twlbbb \scriptfont\mathbbfam=\egtbbb
\scriptscriptfont\mathbbfam=\sixbbb

\makeatother

\newtheorem{theorem}{Theorem}[section]

\newtheorem{lemma}[theorem]{Lemma}
\newfont{\tenbi}{cmbxti10}

\begin{document}

\title{Helicity-type integral invariants for Hamiltonian systems}

\author{Mikhail V. Deryabin \footnote{Mads Clausen Instituttet,
Syddansk Universitet, S{\o}nderborg 6400, Denmark. Email:
mikhail@mci.sdu.dk}}

\date{}

\maketitle

\abstract{ In this note, we consider generalizations of the
asymptotic Hopf invariant, or helicity, for Hamiltonian systems
with one-and-a-half degrees of freedom and symplectic
diffeomorphisms of a two-disk to itself. }

\section{Helicity integral for Hamiltonian systems}

Consider a Hamiltonian system on the extended phase space
$\mathbb{R}^2 \times S^1$ with symplectic coordinates $p,q \in
\mathbb{R}^2$ and time $S^1 = \{t \pmod{2\pi}\}$, given by a
Hamiltonian $H(p,q,t)$:
 \be \label{Ham}
\dot{p} = -H_q, \quad \dot{q} = H_p, \quad \dot{t} = 1.
 \ee
Let $D \in \mathbb{R}^2$ be an invariant disk for this system (which
always exists for integrable systems with bounded energy levels and
near-integrable systems, close to them). An integral over a solid
torus $D \times S^1$:
 \be \label{helicity}
\mathcal{H}(H) = \int_0^{2\pi} \int_{D} \tilde{H}(p,q,t) \, dp
\wedge dq \wedge dt,
 \ee
where the function $\tilde{H} = H + const$ equals zero on the
boundary $\mathbb{T}^2 = \partial D \times S^1$, will be called {\it
helicity} for the Hamiltonian system (\ref{Ham}). It is easy to see
that this definition is consistent -- the value $\mathcal{H}$ does
not depend on the choice of the symplectic coordinates. The term
''helicity'' is chosen, as the integral (\ref{helicity}) equals (up
to a multiplier 2) to the helicity of the Hamiltonian vector field
$\xi = (-H_q, H_p, 1)$ (with invariant measure on the extended phase
space being $\mu = dp \wedge dq \wedge dt$):
 \be \label{hel_proof}
2 \int_{D \times S^1} \tilde{H}(p,q,t) dp \wedge dq \wedge dt =
\int_{D \times S^1} i_\xi \mu \wedge (i_\xi \mu)^{-1},
 \ee
provided the 1-form $(i_\xi \mu)^{-1}$ is chosen such that
 \be \label{condition}
\int_{t \in S^1} (i_\xi \mu)^{-1} = 0.
 \ee
Indeed, as $i_\xi \mu = dp \wedge dq - dH \wedge dt$, then, due to
(\ref{condition}), $(i_\xi \mu)^{-1} = p dq - \tilde{H} dt$.
Relation (\ref{hel_proof}) is now obtained by the integration by
parts, cf.~\cite{MD05}. Notice that, as we have assumed that the
disk $D$ is invariant, from Condition (\ref{condition}) follows that
$\tilde{H} = 0$ on its boundary $\partial D$.

The helicity invariant has the following ergodic interpretation: it
measures asymptotic linking of the trajectories of the Hamiltonian
vector field $(-H_q, H_p, 1)$ in the solid torus $D \times S^1$,
which in turn is itself unlinked (and untwisted), see \cite{MD05}
for details. In particular, the invariant (\ref{helicity}) does not
change after rescaling $H \to \mu H$, $t \to \frac{1}{\mu} t$ -- the
average linking of the trajectories remain the same.

{\bf Remark.} Expression (\ref{helicity}) is still well defined if
the Hamiltonian $H$ is discontinuous at some $t$. One can show that
the ergodic interpretation of (\ref{helicity}) remains the same.

Expression (\ref{helicity}) is exactly the same as the expression
for the Calabi invariant, see \cite{Ca}. The difference is that
above we have not assumed that the gradient of the Hamiltonian $H$
is zero on the boundary torus $\mathbb{T}^2$ (which is an assumption
in the Calabi invariant definition). Of course if we have assumed
that, then (\ref{helicity}) would be the Calabi invariant.

\section{Hamiltonian diffeomorphisms}

Consider a symplectic mapping $h: D \to D$. There exists a
Hamiltonian system (with the Hamiltonian $H$ depending
$2\pi$-periodically on time), which Poincar{\'e} map $g^{2\pi}_H$
coincides with $h$ (actually, there are infinitely many such
systems, and one can show that they have the same smoothness as the
mapping itself, see, e.g., \cite{Tresh}). Notice that for the
calculations below, it does not matter which period to take (due to
the ergodic interpretation of the helicity integral).

The Calabi invariant is an invariant of the symplectic mapping,
identical at the boundary, and it does not depend on a paricular
choice of the underlying Hamiltonian system, see, e.g., \cite{AK}.
On the contrary, if we only demand that the mapping $h$ sends the
boundary $\partial D$ to itself, then the value of the integral
(\ref{helicity}) may depend on a particular choice of the
Hamiltonian flow. However, the values of $\mathcal{H}$ will be
''quantized'' by the square of the sympectic area $S(D)$ of the disk
$D$:

\begin{theorem} \label{Th1}
Let a symplectic mapping $h: D \to D$ be given, which sends a
boundary $\partial D$ to itself. Let $g^t_{H_1}$ and $g^t_{H_2}$ be
two Hamiltonian flows, given by $2\pi$-periodic Hamiltonians $H_1$
and $H_2$, such that $g^{2\pi}_{H_1} = g^{2\pi}_{H_2} = h$. Then,
$\mathcal{H}(H_1) - \mathcal{H}(H_2) = n S(D)^2 / 2$ for some $n \in
\mathbb{Z}$.
\end{theorem}

{\it Proof.} We first prove the following

\begin{lemma} \label{Lem1}
Let $h = id$. Then for any $H$, such that $g^{4\pi}_H = id$,
$\mathcal{H}(H) = n S(D)^2 / 2$, $n \in \mathbb{Z}$.
\end{lemma}

{\it Proof of Lemma \ref{Lem1}.} The mapping $g^{4\pi}_H$ turns the
boundary circle by the angle $2\pi n$, $n \in \mathbb{Z}$ (because
of continuity, no other transformation of the circle than the pure
rotation is allowed). Take a bigger disk $\tilde{D}$, which contains
$D$, and define a mapping $\tilde{h}: \tilde{D} \to \tilde{D}$,
stationary at the boundary $\partial \tilde{D}$, and coinciding with
$h$ in $D$. The Calabi invariant for the mapping $\tilde{h}$ does
not depend on an underlying Hamiltonian flow. If we take the flow
that coincides with $g^t_H$ inside $D$, and that is stationary at
the outer boundary $\partial \tilde{D}$, then, as we tend $\tilde{D}
\to D$, its Calabi invariant will tend to $\mathcal{H}(H)$, as
 $H|_{\partial D}$ will tend to zero as $\tilde{D} \to
D$.

Now, we take the Hamiltonian $\hat{H}$ in $\tilde{D}$, such that in
$D$, $\hat{H} = n I / 2 + c$, where $I$, $\phi$ are the
''action-angle'' variables, such that the disk $D$ is given by
inequality $I \le I_0$, and the disk $\tilde{D}$ is given by $I \le
I_1$, $I_1 > I_0$. We define the constant $c$ later. Let $I_1 = I_0
+ \epsilon$. The function $\hat{H}$ in $\tilde{D} \setminus D$ can
be taken as
 \be \label{H_in_ring}
\hat{H} = \frac{n}{4 \epsilon} (-I + I_0 + \epsilon)^2.
 \ee
Indeed, $\hat{H}|_{\partial \tilde{D}} = 0$, $d\hat{H}|_{\partial
\tilde{D}} = 0$ and
 \bew
\frac{d\hat{H}}{dI} (I_0) = \frac{n}{2 \epsilon} (-I_0 + I_0 +
\epsilon) = \frac{n}{2}.
 \eew
As $\hat{H}(I_0) = {n \epsilon }/{4}$, the constant $c = {n \epsilon
}/{4} - n I_0 / 2$.

The Calabi invariant equals
 \bew
\mathcal{C}(\tilde{h}) = 4\pi \int_{\tilde{D}} \hat{H} dI \wedge
d\phi = 4 \pi^2 n \int_0^{I_0} (I - I_0) dI + O(\epsilon).
 \eew
As we turn $\epsilon \to 0$, we get
 \bew
\mathcal{C}(\tilde{h}) \to -2 \pi^2 n I_0^2 = -\frac{n S(D)^2}{2},
 \eew
where $S(D)$ is the symplectic area of the disk $D$. $\Box$

Now, to prove the theorem, consider a mapping $g^{2\pi}_{H_1} \cdot
g^{2\pi}_{-H_2}$. This mapping is identical. Its asymptotic Hopf
invariant equals sum of asymptotic Hopf invariants for systems with
Hamiltonians $H_1$ and $-H_2$, which, by Lemma \ref{Lem1}, equals $n
S(D)^2 / 2$. $\Box$

{\bf Remark.} The group of Hamiltonian diffeomorphisms of a 2-disk,
identical at the boundary, is known to be contractible, see
\cite{ER}. In particular, the Calabi 1-form, defined as the
right-invariant differential form, coinciding with the Calabi
integral
 \bew
\int_D H dp \wedge dq
 \eew
on the Lie algebra of Hamiltonian diffeomorphisms of disk $D$, is
exact, see \cite{AK}. If we drop the condition that a Hamiltonian
diffeomorphism is identical at the boundary, then from Theorem
\ref{Th1} follows that this form (which is still correctly defined
and closed) is no longer exact, and thus the topology of the group
of Hamiltonian diffeomorphisms becomes more complicated.

The statement converse to Theorem \ref{Th1} is also true. Let a
smooth symplectic mapping $h:D \to D$ and a number $n \in
\mathbb{Z}$ be given.

\begin{theorem} \label{Th2}
There are two smooth Hamiltonian systems, given by Hamiltonians
$H_1$ and $H_2$, $2\pi$-periodic in $t$, such that $h =
g^{2\pi}_{H_1} = g^{2\pi}_{H_2}$, and
 \bew
\mathcal{H}(H_1) - \mathcal{H}(H_2) = n S(D)^2 / 2.
 \eew
\end{theorem}

{\it Proof.} It is well-known that a symplectic mapping can be
inserted into a Hamiltonian flow, see, e.g., \cite{Tresh}. Let this
flow be given by a Hamiltonian $H_1$ (which is zero on the boundary
$\partial D$). Let $I, \phi$ be the ''action-angle'' variables, such
that the disk $D$ is given by $I = I_0$. Take now the following
function:
 \bew
\tilde{H}_2(I,\phi,t) = \Big\{
\begin{array}{cc}
\frac{1}{2} H_1(I,\phi,2t), \quad t \in [0, \pi] \\
2n(I - I_0), \quad t \in [\pi, 2\pi].
\end{array}
 \eew
Obviously, $\mathcal{H}(\tilde{H}_2) - \mathcal{H}(H_1) =
\frac{n}{2} S(D)^2$. Now, we have to construct a smooth Hamiltonian
$H_2$ with the same helicity as $\tilde{H}_2$. We show that this can
be done by time transformation.

Suppose that we want to make $H_2$ differentiable (once, in time).
Consider the following time transformation: $t = \tau - \frac{1}{2}
\sin 2(\tau - \pi)$. The helicity integral becomes:
 \bew
\mathcal{H}(\tilde{H}_2) = \int_0^{2\pi} \int_D (1 - \cos 2(\tau -
\pi)) \tilde{H}_2 (I, \phi, \tau) dI \wedge d\phi \wedge d\tau = \\
= \int_0^{2\pi} \int_D {H}_2 \, dI \wedge d\phi \wedge d\tau =
\mathcal{H}({H}_2),
 \eew
where we denoted $H_2(I, \phi, \tau) = (1 - \cos 2(\tau - \pi))
\tilde{H}_2(I, \phi, \tau)$. Obviously, $H_2 \in C^1$ (being a
$2\pi$-periodic function of $\tau$).

To obtain smoothness of order $k \in \mathbb{N}$, one can construct
time transformation, which is given by $t = \tau + (\mbox{periodic
function of } \tau)$, such that at $\tau = \delta$ and $\tau = \pi +
\delta$, $t = O(\delta^{k+2})$ and $t = \pi + O(\delta^{k+2})$
correspondingly for small $\delta$. For example, for $k = 3$, take
 \bew
t = \tau - \frac{1}{2} \sin 2(\tau - \pi) - \frac{1}{12} \sin^3
2(\tau - \pi). \qquad \Box
 \eew

{\bf Remark.} Theorem \ref{Th2} is also true in the analytic case.
The proof is simple modulus results in \cite{Tresh}: one should for
example use the smoothing technique from \cite{Tresh} to get an
analytic system from the smooth one, obtained during the proof of
Theorem \ref{Th2}. Notice also that if two Hamiltonian functions are
close to each other (in an appropriate topology) and they define the
same mapping, then the helicity will also be the same -- by
continuity of the helicity functional and by Theorem \ref{Th1}.

From Theorem \ref{Th2} follows that a symplectic mapping can be
inserted into a Hamiltonian flow with any given helicity level. In
particular, there exists the flow with smallest (in the absolute
value) helicity, which is a functional of the mapping only. This
smallest helicity will be called the {\it generalized Calabi
invariant}.

We do not write down explicitly the generalized Calabi invariant as
the functional of a symplectomorphism here. We only note that it can
be written as a sum of two parts: an integral over the disk, and an
integral over the boundary circle. One way to see that is to take
the limiting procedure, as we did in the proof of Lemma \ref{Lem1}.

\end{document}